\documentclass[11pt]{amsart}
\usepackage{amssymb, amsmath, amscd, dcpic, pictexwd, mathtools}
\usepackage{hyperref}

\usepackage{amsopn}

\usepackage{tikz-cd}

\usepackage{leftidx}

\usepackage{xcolor}

\newtheorem{thm}{Theorem}[section]

\newtheorem{rem}[thm]{Remark}

\numberwithin{equation}{section}

\newcommand{\Z}{\mathbb{Z}}
\newcommand{\Q}{\mathbb{Q}}

\author{An Huang, Christian Jepsen}

\title{One-loop $p$-adic string theory \\ 
and the N\'eron local height function}
\begin{document}
\maketitle
\begin{abstract}
The $p$-adic string worldsheet action on the quotient of the Bruhat-Tits tree of $PGL(2,\mathbb{Q}_p)$ by a genus 1 Schottky group has a dual description on the asymptotic boundary, the Tate curve $\mathbb{Q}_p^\ast/q^\mathbb{Z}$ \cite{HJ2024}. We show that the two point function of the dual action coincides with the N\'eron-Tate local height function of the Tate curve. 
\end{abstract}

\begin{center}
\hspace{1mm}
\\
\emph{Dedicated to the memory of Jo\"el Bella\"iche}
\\[6mm]
\end{center}

\section{Introduction}

The recent work \cite{HRSW2025} showed that the N\'eron-Tate local height function on the Tate curve is equal to a large $p$ limit of the two-point function of an action defined on the Tate curve when the $j$-invariant has odd valuation. The action defined in \cite{HRSW2025} was inspired by $p$-adic string theory. The aim of the present paper is to refine this result by clarifying the relation between the height function and $p$-adic strings for any $p$.

The study of $p$-adic strings first derived motivation from the observation of Freund and Olson \cite{FO1987} in 1987 that scattering amplitudes in this formalism factorize and exhibit dual resonance and crossing symmetry in accordance with physical expectations. Subsequently Brekke, Freund, Olson, and Witten \cite{Brekke1988} were able to explicitly evaluate all tree-level $N$-point amplitudes and write down an effective spacetime action that was later used by Ghoshal and Sen \cite{Ghoshal2000} and by Gerasimov and Shatashvili \cite{Gerasimov2000} to argue that the $p$-adic formalism sheds light on the condensation of tachyons and their relation to $D$-branes.

The tree-level worldsheet action of $p$-adic string theory was formulated in 1989 by Zabrodin \cite{Z1989}, who demonstrated that from a theory of a graph Laplacian acting on the Bruhat-Tits tree $T_p$ of the group $PGL(2,\mathbb{Q}_p)$ it is possible to derive the Freund-Olson amplitudes as well the equivalent action that had earlier been constructed by Zhang \cite{Zhang1988} and Spokoiny \cite{Spokoiny1988}, and which is defined over the boundary of the open string worldsheet. The dual action in this latter boundary formulation is given, up to a choice of normalization constant, as
\begin{equation}
S_1=\int_{\mathbb{Q}_p}\phi(x)D_1\phi(x)\,dx\,,
\end{equation}
where $D_1$ is the regularized Vladimirov derivative defined by the following integral operator on the space of compactly supported locally constant functions on $\mathbb{Q}_p$:
\begin{equation}
D_1\phi(x):=\int_{\mathbb{Q}_p}\frac{\phi(z)-\phi(x)}{|z-x|^2}\,dz\,.
\end{equation}
$D_1$ is a pseudo-differential operator whose symbol is given by the $p$-adic norm function. $D_1$ and a family of its deformations have been studied in numerous papers, with ideas/results partially summarized in \cite{huang2021greens}\cite{HSZ2022}\cite{HRSW2025}. In particular, the operator is closely related to dimensional regularization, Hecke-L functions, Tate's thesis, quadratic reciprocity, non-Archimedean harmonic analysis, and Verma modules. 

\begin{figure}
\centering
$
\begin{matrix}
\text{
\includegraphics[width=0.4\columnwidth]{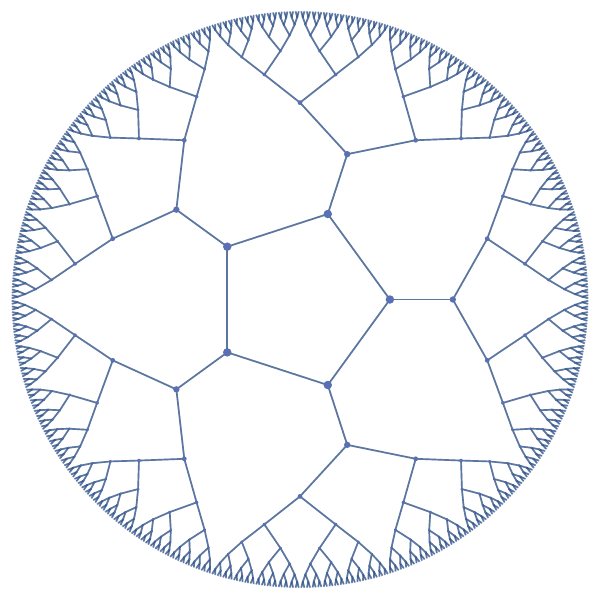}}
\end{matrix}
$
\caption{The $p$-adic string worldsheet at one loop: the tree quotient $T_p/\Gamma$. For this example $p=2$ and $m=5$. \label{Fig:TreeQuotient} 
} 
\end{figure} 

As a next step, one is interested in $p$-adic strings in genus one. These were first explored in \cite{CMZ1989}, where the one-loop worldsheet is a quotient of the tree: $T_p/\Gamma$, with $\Gamma$ being the discrete subgroup of $PGL(2,\mathbb{Q}_p)$ generated by $\begin{bmatrix} q & 0\\0 & 1\end{bmatrix}$, and $|q|<1$. As Zabrodin's action only depends on $q$ through $|q|$, we assume $q=p^m$. See Figure~\ref{Fig:TreeQuotient} for an example of this tree graph quotient. In \cite{HJ2024}, a dual action on $\Q_p^*/q^{\Z}$, the asymptotic boundary of $T_p/\Gamma$, has been worked out: 
\begin{equation}\label{dual action}
S=\int_E\phi(x)D\phi(x)d^*x
\end{equation}
where the operator $D$ appearing in the action is defined as
\begin{equation}\label{Ddef}
\hspace{25mm}
D\phi(x):=-c_p\int_E H(z,x)(\phi(z)-\phi(x))d^*z\,,
\hspace{10mm}
c_p:=\frac{p(1-p^{-1})^2}{1-p^{-2}} \,,
\end{equation}
with the function $H(z,x)$ defined as
\begin{align}
\label{Hdef}
H(z,x) := \frac{|x|\,|z|}{|x-z|^2}
+\frac{1}{p^m-1}\Big(\frac{|x|}{|z|}+\frac{|z|}{|x|}\Big)\,,
\end{align}
where we normalize the additive Haar measure $dx$ so that the measure of $\mathbb{Z}_p$ is $1$. The multiplicative Haar measure is normalized by $d^*x=\frac{dx}{|x|}$. $E$ is a fundamental domain of the $\mathbb{Q}_p$ points of the Tate curve:
\begin{align}
E:=\cup_{i=0}^{m-1}p^i\mathbb{Z}_p^\times\,,
\end{align}
where $\mathbb{Z}_p^\times$ are the $p$-adic units (also sometimes denoted $\mathbb{U}_p$).

It is straightforward to verify that $D$ is self-adjoint under the integral paring 
\begin{align}
<f(x),g(x)>:=\int_Ef(x)g(x)d^*x\, 
\end{align}
and that $D$ is a positive semi-definite operator whose kernel is given by constant functions. We shall see that $D$ is also a pseudo-differential operator, with a spectral gap and Weyl asymptotics in line with what one would expect from a Laplacian operator in two dimensions. We explicitly compute the spectrum of $D$ and its regularized determinant. 

The main point of this paper is to show that up to an undetermined additive constant, the Green's function $G(x,y)$ for $D$, i.e. the two point function of \eqref{dual action}, coincides with the N\'eron-Tate local height function $h(x)$, given below in equation~\eqref{h}. More precisely, we show that up to an additive constant, there is a unique symmetric Green's function $G(x,y)$, and $G(x,y)=h(x/y)$ when $v(x)\geq v(y)$. Thus, in terms of the relation between $p$-adic strings and the height function, the result of this paper refines that of \cite{HRSW2025}. In consequence, the local height function as such a Green's function represents a basic building block for computing $p$-adic string amplitudes on the Tate curve. On the other hand, the result implies a role of physics, in particular $p$-adic strings, in the foundations of arithmetic geometry. 

\begin{rem}
It is straightforward to extend the result to a finite extension of $\Q_p^*$. We only state the result for $\Q_p^*$.
\end{rem}

\section{Properties of the dual action}
\label{sec:Properties}
The dual action derived in \cite{HJ2024} was not originally written as \eqref{dual action} with the weight function $H(z,x)$ given as in \eqref{Hdef}. Rather, letting $v(x)$ denote the valuation of $x\in E$ so that $|x|=p^{-v(x)}$, the weight function was written in the equivalent form
\begin{align}
\label{Hcases}
H(z,x)
=\,&
\begin{cases}
\displaystyle
\frac{|x|\,|z|}{|z-x|^2}+\frac{2}{p^m-1}
&\text{for }v(z)=v(x)\,,
\\[8mm]
\displaystyle
\frac{p^{m-u}+p^u}{p^m-1}
&\text{for }
0<|v(z)-v(x)|=u<m\,,
\end{cases}
\end{align}
which is sometimes useful for practical computations. The advantage of the form \eqref{Hdef} is that it enables us more readily to demonstrate some of the symmetry properties of the action. Specifically, the action is invariant under dilatations and inversions, as we will presently explain in detail.

By a dilatation, we mean the operation of replacing $\phi$ with a function $\widetilde{\phi}$ obtained from $\phi$ by multiplication of the argument by any fixed element $\lambda \in E$: $\widetilde{\phi}(x) = \phi(\lambda x)$. Since the action is given in terms of integration over the multiplicative Haar measure, which is invariant under rescalings $x\rightarrow \lambda x$, the dual action $S$ is transformed by a dilatation into the following (times a constant prefactor):
\begin{align}
\label{actionEE}
\int_{E\times E}d^\ast x\,d^\ast z
\,\phi(\lambda x)\,H(z,x) \Big(\phi(\lambda x)-\phi(\lambda z)\Big)
=
\int_{E\times E}d^\ast x\,d^\ast z
\,\phi(x)\,H(\frac{z}{\lambda},\frac{x}{\lambda}) \Big(\phi(x)-\phi(z)\Big)\,.
\end{align}
From this equation, we see that dilatational invariance of the action will immediately follow if we can establish that $H(z,x)$ equals $H(\frac{z}{\lambda},\frac{x}{\lambda})$ for all $\lambda,z,x\in E$. At a first glance, this might seem immediately obvious from \eqref{Hdef} since each term is given by a ratio of norms. But the identity is slightly more subtle and in fact the precise form of the prefactor $(p^m-1)^{-1}$ in \eqref{Hdef} is crucial. The issue is that the ratios $z':=z/\lambda$ and $x':=x/\lambda$ in \eqref{actionEE} designate elements in $E$, which means that in some instances one must act with the quotienting group $\Gamma$ to ensure that a ratio remains inside the fundamental domain. Phrased differently, since each element in $E$ has a norm between 1 and $p^{-m}$, the norm of a ratio need not equal the ratio of norms. Concretely, suppose $|z|>|x|$, which implies that $|z-x|=|z|$. In this case there exists $\lambda\in E$ such that $|x'|=\frac{|x|}{|\lambda|}$ while $|z'|=\frac{|z|}{|\lambda|}p^{-m}$ (for example, one could pick $\lambda=x$). In this case, we have that $|z'|<|x'|$, which implies that $|z'-x'|=|x'|$, and so it follows that
\begin{align}
H(z',x') =\,&
\frac{|x'|\,|z'|}{|x'|^2}
+\frac{1}{p^m-1}\Big(\frac{|x'|}{|z'|}+\frac{|z'|}{|x'|}\Big)
\nonumber\\[1mm]
=\,&\frac{|z|}{|x|}p^{-m}
+\frac{1}{p^m-1}\Big(\frac{|x|}{|z|p^{-m}}+\frac{|z|p^{-m}}{|x|}\Big)
\nonumber\\[-2.3mm]
\\[-2.3mm]
=\,&
\frac{|x|}{|z|}
+\frac{1}{p^m-1}\Big(\frac{|x|}{|z|}+\frac{|z|}{|x|}\Big)
\nonumber\\[1mm] \nonumber
=\,&
\frac{|x|\,|z|}{|x-z|^2}
+\frac{1}{p^m-1}\Big(\frac{|x|}{|z|}+\frac{|z|}{|x|}\Big)
=H(z,x)\,.
\end{align}
The conclusion is that the identity $H(z,x)=H(\frac{z}{\lambda},\frac{x}{\lambda})$ really is true for all $x,y,\lambda \in E$, and so the dilatational invariance of the action follows.

By an inversion, we mean the operation of replacing $\phi$ with a function $\widehat{\phi}$ obtained from $\phi$ by taking the reciprocal of the argument: $\widehat{\phi}(x) = \phi(\frac{1}{x})$. Since the multiplicative Haar measure is invariant under reciprocation, the invariance of the action is established by demonstrating that $H(z,x)=H(\frac{1}{z},\frac{1}{x})$ for all $z,x\in E$. If $|x|=|z|$ or if $|x|$ and $|z|$ are both less than one, inversional invariance is true separately for $\frac{|x|\,|z|}{|x-z|^2}$ and $\frac{|x|}{|z|}+\frac{|z|}{|x|}$. For the remaining type of situation when, say, $|x|=1$ and $|z|<1$ (which implies that $|x-z|=1$ and $|\frac{1}{z}|=\frac{1}{|z|}p^{-m}$) the precise form of the coefficient $(p^m-1)^{-1}$ again becomes crucial:
\begin{align}
\nonumber
H(\frac{1}{z},\frac{1}{x})=\,&
p^{-m}|z|^{-1}
+\frac{1}{p^m-1}
\Big(|z|p^m+\frac{1}{|z|p^m}\Big)
\\
=\,&
|z|+\frac{1}{p^{m}-1}
\Big(\frac{1}{|z|}+|z|\Big)
\\
=\,&
\frac{|x|\,|z|}{|x-z|^2}
+\frac{1}{p^m-1}\Big(\frac{|x|}{|z|}+\frac{|z|}{|x|}\Big)
=H(z,x)\,.
\nonumber
\end{align}
We conclude that $H(z,x)=H(\frac{1}{z},\frac{1}{x})$ for every $z,x\in E$, which implies that the action is invariant under inversions.

\begin{rem}
Note that the precise form of the function $H(x,z)$ in the dual action, which was determined in \cite{HJ2024} by a direction computation, could alternatively have been derived solely from first principle symmetry considerations: first, the leading (in the limit as $z$ approaches $x$) term $|x||z|/|x-z|^2$ is determined by the conformal dimension and the rotation symmetry; then the second term $\frac{1}{p^m-1}\Big(\frac{|x|}{|z|}+\frac{|z|}{|x|}\Big)$ has to be added to account for how the symmetry interacts with the fundamental domain $E$, as explained.
\end{rem}

\begin{rem}
\label{remark:Covariance}
While the properties $H(x,y)=H(\lambda x,\lambda y)$ and $H(x,z)=H(\frac{1}{x},\frac{1}{z})$ imply dilational and inversional invariance of the action, for the operator $D$, they imply inversional and dilational \emph{covariance}. That is to say, $D[\phi(\lambda \cdot)](x)=D[\phi(\cdot)](\lambda x)$ and $D[\phi(\frac{1}{\cdot})](x)=D[\phi(\cdot)](\frac{1}{x})$.  
\end{rem}
 
\section{Green's function and the height function}
\label{sec:Green}
For $x \in E$, let $v_x=v(x)$ and let 
\begin{align}\label{h}
h(x)=v(x-1)+\frac{v_x(v_x-m)}{2m}+\frac{m}{12}
\end{align}
be the N\'eron-Tate local height function on the Tate curve \cite{Silverman}. In the following, we shall show that $G(x,1):=h(x)$ is a Green's function of $D$ when the second argument is $1$. To show this, we will verify by direction computation that 
\begin{align}
\label{DHeq}
Dh(x)=-\frac{1}{\text{Vol}}=-\frac{p}{m(p-1)}\hspace{3mm}\text{for}\hspace{3mm}x\neq 1\,,
\end{align}
where $\text{Vol}$ is the Volume of the Tate curve with respect to the multiplicative Haar measure. Since the kernel of $D$ consists of constant functions, such a verification suffices to imply that $h(x)$ is a Green's function of $D$, by the same argument as in \cite{HRSW2025}. 

To compute $Dh(x)$ according to the definition in \eqref{Ddef}, we must perform an integral over all $z\in E$. We can partition this integral according to the valuation of $z$ as
\begin{align}
\hspace{10mm}
-\frac{1}{c_p}Dh(x)
=\sum_{v_z=0}^{m-1}
I_{v_z}
\hspace{5mm}
\text{with}
\hspace{5mm}
I_{v_z}:=
\int_{p^{v_z}\mathbb{Z}_p^\times} d^\ast z\,
 H(z,x)
 \big(
h(z)-h(x)
 \big)\,.
\end{align} To evaluate the different pieces in this decomposition, we find it convenient to consider in turn the two cases: $v_x =0$ and $v_x \neq 0$, assuming in the former case that $x \neq 1$.
\\[2mm]
\emph{Case 1: $v_x =0$.} When $v_z =0$, the top situation in \eqref{Hcases} applies. Then, letting $l=v(x-1)$ and $j=v(z-1)$,
\begin{align}
\nonumber
I_0=\,
& 
\int_{\mathbb{Z}_p^\times} d^\ast z\,
\Big(
\frac{1}{|z-x|^2}+\frac{2}{p^{m}-1}
\Big)
  \Big(
j-l
 \Big)
\\ 
=\,&
\frac{p-2}{p}
\Big(
1+\frac{2}{p^{m}-1}
\Big)
(-l)
+
\frac{p-1}{p}\sum_{j=1}^{l-1}
p^{-j}
\Big(
p^{2j}+\frac{2}{p^{m}-1}
\Big)
  \Big(
j-l
 \Big)
\\ \nonumber
&\hspace{40.2mm}
+
\frac{p-1}{p}\sum_{j=l+1}^\infty
p^{-j}\Big(
p^{2l}+\frac{2}{p^{m}-1}
\Big)
  \Big(
j-l
 \Big)\,.
\end{align}
Meanwhile, when $z$ has a non-zero valuation, $v(z-1)=0$ and the integral pieces evaluate as
\begin{align}
I_{v_z>0} =\,& \int_{p^{v_z}\mathbb{Z}_p^\times} d^\ast z\,
 \frac{p^{m-v_z}+p^{v_z}}{p^m-1}
 \Big(
\frac{v_z(v_z-m)}{2m}-l
 \Big)
 \nonumber \\[-2.5mm]
 \\[-2.5mm] \nonumber
=\,&
\frac{p-1}{p}
 \frac{p^{m-v_z}+p^{v_z}}{p^m-1}
 \Big(
\frac{v_z(v_z-m)}{2m}-l
 \Big)\,.
\end{align}
Adding together the different contributions to the integral for $Dh(x)$ and carrying out the finite and infinite summations, the dependency on $l$ can be seen to cancel out and the advertised result follows,
\begin{align}
-\frac{1}{c_p}Dh(x)
=I_0 + \sum_{v_z=1}^{m-1}I_{v_z}
=\frac{p+1}{m(p-1)^2}
=\frac{1}{c_p\text{Vol}}\,.
\label{case1ans}
\end{align}
\emph{Case 2: $0 < v_x < m$.} In this case $v(x-1)=0$. We first consider the sub-case $v_z =0$ so that $v_z \neq v_x$ and the bottom situation in \eqref{Hcases} applies. Denoting $v(z-1)=j$ again, we have that
\begin{align}
I_0 =\,&
\int_{\mathbb{Z}_p^\times} d^\ast u_z\,
\frac{p^{m-v_x}+p^{v_x}}{p^m-1}
  \Big(
j-\frac{v_x(v_x-m)}{2m}
 \Big)
\nonumber  \\
 =\,&
\frac{p^{m-v_x}+p^{v_x}}{p^m-1}
  \Big(
\frac{p-1}{p}\sum_{j=1}^\infty j p^{-j}
-\frac{p-1}{p}\frac{v_x(v_x-m)}{2m}
 \Big)
   \\ \nonumber
 =\,&
\frac{p^{m-v_x}+p^{v_x}}{p^m-1}
  \Big(
\frac{1}{p-1}
-\frac{p-1}{p}\frac{v_x(v_x-m)}{2m}
 \Big)\,.
\end{align}
Meanwhile, the integral piece where $x$ and $z$ have the same valuation vanishes, $I_{v_z=v_x}=0$, since in this case $h(z)-h(x)=0$. For the remaining integral pieces, with $0<v_z<v_x$ and $v_x<v_z$, we have that
\begin{align}
I_{v_z<v_x}=\,&\frac{p-1}{p}\,
 \frac{p^{m-v_x+v_z}+p^{v_x-v_z}}{p^m-1}
 \Big(
\frac{v_z(v_z-m)}{2m}-\frac{v_x(v_x-m)}{2m}
 \Big)\,,
 \nonumber \\[-2.5mm]
 \\[-2.5mm] \nonumber
I_{v_z>v_x}=\,&
\frac{p-1}{p}\,
 \frac{p^{m-v_z+v_x}+p^{v_z-v_x}}{p^m-1}
 \Big(
\frac{v_z(v_z-m)}{2m}-\frac{v_x(v_x-m)}{2m}
 \Big)\,.
\end{align}
Adding together all the integral pieces, we find that
\begin{align}
-\frac{1}{c_p}Dh(x)
=I_0 
+ \sum_{v_z=1}^{v_x-1}I_{v_z}
+ \sum_{v_z=v_x+1}^{m-1}I_{v_z}
=\frac{p+1}{m(p-1)^2}
=\frac{1}{c_p\text{Vol}}\,.
\end{align}
And so we conclude that in both cases $Dh(x) = -1/\text{Vol}$.

We have now established that $G(x,1)=h(x)$ is a Green's function of $D$. Next, we extend the result to a Green's function $G(x,y)$ where the second argument is arbitrary. This extension follows immediately from the property, established in Section~\ref{sec:Properties} and commented upon in Remark~\ref{remark:Covariance}, that $D$ transforms covariantly under dilatations. By this property, we conclude that for any $y \in E$, the function $G(x,y):=h(\frac{x}{y})$ is a Green's function of $D$, ie. $DG(x,y) = -1/\text{Vol}$ for any $x\neq y$. Note that $\frac{x}{y}$ is understood to represent a unique point in the fundamental domain $E$, as before. And from the inversional covariance of $D$, also established in Section~\ref{sec:Properties}, we conclude that the function $h(\frac{y}{x})$ likewise is a Green's function for $D$.

So by general arguments as in \cite{HRSW2025}, $h(\frac{x}{y})$ is the unique symmetric Green's function for $D$ up to an addition of a constant: i.e. for any $y$, $Dh(\frac{x}{y})+\frac{1}{\text{Vol}}$ is a distribution on the space of locally constant functions supported at $y$, thus it has to be a multiple of $\delta_{x,y}$. The multiple has to be 1 as the integral of $(\delta_{x,y}-\frac{1}{\text{Vol}})(1)$ over $E$ is zero, as 1 is in the kernel of the self-adjoint operator $D$. We therefore have the following:

\begin{thm}
    $G(x,y):=h(\frac{x}{y})$ when $m-1\geq v(x)\geq v(y)\geq 0$, extended by symmetry $G(x,y)=G(y,x)$ to the case when $v(x)< v(y)$, is the unique symmetric Green's function up to an additive constant, for the operator $D$, i.e. $DG(x,y)=\delta_{x,y}-\frac{1}{\text{Vol}}$.
\end{thm}

\begin{rem}
For any point $x$ on the generic fiber of the minimal regular model of the Tate curve, $x$ extends uniquely to a horizontal section $\tilde{x}$ by the valuative criterion for properness. Up to an additive constant, we thus have the interpretation 
\begin{align}
G(x,y)=-\frac{1}{2}\left<\tilde{x}-\tilde{y}+V,\tilde{x}-\tilde{y}+V\right>\,,
\end{align}
where $\left<\,,\right>$ is the intersection pairing on the minimal regular model, and $V$ is a vertical divisor supported at the special fiber, so that the ``corrected divisor" $\tilde{x}-\tilde{y}+V$ is orthogonal to vertical divisors. When $y=1$, this reduces to the definition of the local height, in terms of the intersection pairing.
\end{rem}

\begin{rem}
When $m=1$, the operator $\frac{1}{c_p}\,D$ of the present paper is given by a kernel that is the same as that in \cite{HJ2025} up to a constant shift. Thus our $D$ has the same eigenfunctions and multiplicities as in there, and the spectrum of our $D$ is $\lambda_n=c_p\,(\lambda_n'+\frac{2}{p^m-1}\mu^*(\mathbb{Z}_p^\times))$, where $n$ denotes the conductor of the character, $\lambda_n'=(p+1)p^{n-2}-\frac{2}{p}$ denotes the corresponding eigenvalue as in \cite{HJ2025}, and $\mu^*(\mathbb{Z}_p^\times)=\frac{p-1}{p}$ is the volume of $\mathbb{Z}_p^\times$. Plugging in $m=1$, we get $\lambda_n=(p-1)p^{n-1}$. Therefore, when $m=1$, $\frac{1}{p}D$ agrees exactly with the Dirichlet-to-Neumann boundary operator defined as in \cite{HJ2025}.
\end{rem}

\section{Spectrum of $D$}
The spectrum of the operator $D$ is given by the set of multiplicative characters $\pi$ on $E$. This follows from the multiplicative invariance established in Section~\ref{sec:Properties} that $H(x,z)=H(\lambda\,x,\lambda\,z)$ for all $\lambda,x,z\in E$. For, performing a change of variable $z \rightarrow zx$ and using the property $\pi(zx) = \pi(z)\,\pi(x)$, we have that
\begin{align}
\nonumber D\pi(x)
=\,&
-c_p\int_E d^\ast z\,
H(z,x) \Big(\pi(z)-\pi(x)\Big)
\\
=\,&
-c_p\int_E d^\ast z\,
H(zx,x) \Big(\pi(zx)-\pi(x)\Big)
\\ \nonumber 
=\,&
-\pi(x)\,c_p\int_E d^\ast z\,
H(z,1) \Big(\pi(z)-1\Big)\,,
\end{align}
from which we see that any multiplicative character $\pi$ is an eigenfunction of $D$ with eigenvalue
\begin{align}
\label{lambdaPi}
\lambda_\pi =
-c_p \int_E d^\ast z\,
H(z,1) \Big(\pi(z)-1\Big)\,.
\end{align}
Any element $z \in E$ can be uniquely written as $z = p^{v_z}\,\hat{z}$ with $v_z \in \mathbb{Z}/(m\mathbb{Z})$ and $\hat{z}\in \mathbb{Z}_p^\times$, and the characters over $E$ decompose into a product of characters $\zeta$ over $\mathbb{Z}/(m\mathbb{Z})$ and characters $\hat{\pi}$ over $\mathbb{Z}_p^\times$. Therefore, eigenfunctions of $D$ are given by products of characters along the radial direction and along the circle direction: $\pi(z) = \zeta(v_z)\,\hat{\pi}(\hat{z})$. This is in analogy with the Archimedean case. On the other hand --- as we will show presently --- whenever the radial direction character $\hat{\pi}$ is nontrivial, the eigenvalue only depends on the radial direction character, and it depends only on its conductor. This is different from the Archimedean case.

The fact that a non-trivial character $\hat{\pi}$ implies that $\lambda_\pi$ is independent from the character $\zeta$, can be seen from \eqref{lambdaPi} by decomposing the integral over $E$ into integrals over $p^{v_z}\mathbb{Z}_p^\times$ with $v_z\in\{0,...,m-1\}$. For whenever $v_z \neq 0$, the value of $H(z,1)$ depends only on $|z|$ and not on $\hat{z}$. Therefore, using equation \eqref{Hcases} and the decomposition $\pi(z)=\zeta(v_z)\,\hat{\pi}(\hat{z})$, we have that, for $v_z\neq 0$,
\begin{align}
\int_{p^{v_z}\mathbb{Z}_p^\times}
d^\ast z\,
H(z,1)\,\pi(z)
=
\frac{p^{m-v_z}+p^{v_z}}{p^m-1}\,
\zeta(v_z)
\int_{\mathbb{Z}_p^\times}
d^\ast \hat{z}\,\hat{\pi}(\hat{z})
=0\,.
\end{align}

Since $\zeta(\overline{0})=1$ for any character $\zeta$ on $\mathbb{Z}/(m\mathbb{Z})$, the contribution to the integral \eqref{lambdaPi} for $v_z = 0$ is given by
\begin{align}
-c_p\int_{\mathbb{Z}_p^\times}
d^\ast z\,H(z,1)\Big(\pi(z)-1\Big)=\,&
-c_p\int_{\mathbb{Z}_p^\times}
d^\ast z\,\Big(\frac{1}{|z-1|^2}+\frac{2}{p^m-1}\Big)\Big(\hat{\pi}(z)-1\Big)
\nonumber \\
=\,&
c_p\int_{\mathbb{Z}_p^\times}
d^\ast z\,\frac{1-\hat{\pi}(z)}{|z-1|^2}+c_p\, \mu^*(\mathbb{Z}_p^\times)\,\frac{2}{p^m-1}
\,.
\end{align}
In \cite{HJ2025}, it was shown that the integral $\int_{\mathbb{Z}_p^\times}
d^\ast z\,\frac{1-\hat{\pi}(z)}{|z-1|^2}$ depends only on the conductor $n$ of $\hat{\pi}$ and is equal to $(p+1)p^{n-2}-2/p$. We conclude that for any character $\pi$ whose radial part has conductor $n>0$, the eigenvalue is given by
\begin{align}
\label{radialPi}
\lambda_n
=\,& c_p\bigg(
\frac{p+1}{p^2}p^n-\frac{2}{p} + \frac{p-1}{p}\,\frac{2}{p^m-1} +
\frac{p-1}{p}
\sum_{v_z=1}^{m-1}
\frac{p^{m-v_z}+p^{v_z}}{p^m-1}
\bigg)
=(p-1)p^{n-1}\,.
\end{align}
We observe that the $m$ dependencies cancel in the calculation. The degeneracies of the eigenvalues \eqref{radialPi}, determined by multiplying the number of characters of conductor $n$ with the angular degeneracy due to the $\zeta$ characters, equal $m(p-2)$ for $n=1$ and $m(p-1)^2p^{n-2}$ for $n>1$.

Next let us determine the eigenvalues of the characters with trivial radial character but non-trivial angular character. Let $\zeta$ be a character of the finite group $\mathbb{Z}/(m\mathbb{Z})$ such that $\zeta(\bar{1})=\omega$, where $\omega$ is an $m$-th root of unity. Then a direct computation shows that $\zeta$ is an eigenfunction for $D$ with eigenvalue
\begin{align}
\nonumber
\lambda_{\omega}&=-c_p\,\sum_{v_z=1}^{m-1}\frac{p^{v_z}+p^{m-v_z}}{p^m-1}(\omega^{v_z}-1)\,\mu^*(\mathbb{Z}_p^\times)\\
\label{lamOmega}
&=-c_p\,\frac{p-1}{p}\Big(\frac{1}{p\,\omega-1}+\frac{1}{p\,\omega^{-1}-1}-\frac{2}{p-1}\Big)\\
&=\frac{p(p-1)(2-\omega-\overline{\omega})}{p^2-p(\omega+\overline{\omega})+1}\,.
\nonumber
\end{align}
So the multiplicity of $\lambda_{\omega}$ is 2 if $\omega$ is not real, and the multiplicity is 1 otherwise, which only happens when $m$ is even and $\omega=-1$. 

From the foregoing discussion, it follows that $D$ is diagonalized by continuous characters of the Tate curve, i.e. $D$ is a pseudo differential operator on the Tate curve. 

$\lambda_{\omega}$ is a decreasing function in $\omega+\overline{\omega}$, and for all primes $p$ and for all $m$ the smallest angular eigenvalue is smaller than the smallest radial eigenvalue of $p-1$. So the spectral gap of $D$ is 
\begin{align}
\lambda_{\omega=e^{2\pi i/m}}=\frac{p(p-1)(2-2\cos\big(\frac{2\pi}{m})\big)}{p^2-2p\cos(\frac{2\pi}{m})+1}\,,
\end{align}
which to leading order in $1/m$ equals $\frac{p}{p-1}(\frac{2\pi}{m})^2$.

\begin{rem}
The following Weyl's law is observed for $D$:

Take a big $\lambda>0$ that is much bigger than $p$. Let $M$ be the biggest integer such that $\lambda_{n=M}\leq \lambda$. Then the number of eigenvalues less than or equal to $\lambda$ is
\begin{equation}
N(\lambda)=m(p-1)^2\sum_{i=0}^{M-2}p^i+m(p-2)+m=m\lambda_M\,.
\end{equation}

So, the spectrum of $D$ does behave as what one would expect from the Laplacian of a two-dimensional domain, even though the eigenfunctions degenerate in a different way compared to the Archimedean case.
\end{rem}

\section{Determinant of $D$}

In a physical system whose action is given in terms of a differential operator
$\mathcal{D}$, the vacuum energy energy is computed via the operator determinant, given by the product over non-zero eigenvalues $\lambda$ of $\mathcal{D}$ whenever this product converges, 
\begin{align}
\det \mathcal{D}
= \prod_{\lambda} \lambda\,,
\end{align}
where the product is counted with multiplicity. When the product does not converge, an operator determinant can instead be defined via analytic continuation of a zeta function $\zeta_\mathcal{D}$ associated to $\mathcal{D}$,
\begin{align}
\det \mathcal{D}
=e^{-\frac{d}{ds}\zeta_\mathcal{D}(s)|_{s=0}}
\,,
\hspace{20mm}
\zeta_\mathcal{D}(s):= 
\sum_{\lambda} \lambda^{-s}\,.
\end{align}
In the case of the 1-loop vacuum energy of the $p$-adic string, it turns out it is possible to explicitly evaluate the zeta-regulated determinant of $D$. 

We first note that the angular eigenvalues \eqref{lamOmega} are finite in number, being associated to $m-2$ different $m$-th roots of unity. For this reason, these eigenvalues make a finite contribution to $\det D$, which can be determined independently of the zeta function regularization. Moreover, this contribution can be evaluated in closed form:
\begin{align}
\label{angularProd}
\prod_{\omega}\lambda_\omega
=
\prod_{\ell=1}^{m-1}\frac{2p(p-1)(1-\cos(\frac{2\pi \ell}{m}))}{p^2-2p\cos(\frac{2\pi \ell}{m})+1}
=\frac{m^2 (p - 1)^{m + 1} p^{m - 1}}{(p^m - 1)^2}\,.
\end{align}
Next we turn to the radial eigenvalues \eqref{radialPi} and define a zeta function associated to the multiplicative characters with non-trivial radial character $\hat{\pi}$:
\begin{align}
\nonumber
\zeta_\pi(s)=\,&
\sum_{\pi \text{ with non-trivial }\hat{\pi}}\lambda_n^{-s}
\\ 
\label{zetaPi}
=\,&
m(p-2)\,(p-1)^{-s}
+\sum_{n=2}^\infty m(p-1)^2p^{n-2}\,\Big((p-1)p^{n-1}\Big)^{-s}
\\
=\,&
\frac{m(p^{s+1}-2p^s+1)}{(p^s-p)(p-1)^s}\,.
\nonumber
\end{align}
While the sum over $n$ in \eqref{zetaPi} converges only for $s>1$, the resummed answer straightforwardly admits analytic continuation in $s$, from which we compute the contribution to $\det D$,
\begin{align}
\label{radialContribution}
e^{-\zeta_\pi'(s)|_{s=0}}=\left(\frac{p}{p-1}\right)^m\,.
\end{align}
Multiplying together the contribution \eqref{angularProd} from the purely angular characters with the contribution \eqref{radialContribution} from characters with a radial dependence, we arrive at the zeta-regulated determinant
\begin{align}
\label{detD}
\det D = m^2\frac{1-p^{-1}}{(1-p^{-m})^2}\,.
\end{align}
It remains an open problem to ascertain if there exists a sensible method of combining together the partition functions of distinct graphs into a single partition function for fluctuating graph geometries. If such a model of quantum gravity exists, the determinant \eqref{detD} will enter into a suitably weighted sum over $m$ for the computation of the full one-loop gravitational partition function.

\section{The N\'eron local height function and holography}

From the definition of the operator $D$ in \eqref{Ddef}, it is clear that this operator is insensitive to any constant shift in the function it acts on, and so the stipulation that a function $h(x)$ is a Green's function of $D$ does not permit a determination of the constant term $\frac{m}{12}$ in \eqref{h}. However, there is a specific sense in which this precise term is encoded in the holographic $p$-adic AdS/CFT correspondence, which was first advanced in \cite{Gubser2017} and \cite{HMST}, and which gave rise to the dual action \eqref{dual action}. 

The action \eqref{dual action} furnishes the boundary dual to the bulk theory on $T_p/\Gamma$ containing a kinetic term given by the graph Laplacian. The dual action was derived in \cite{HJ2024}, which contained a generalization of the dual action to the case when the bulk action contains a non-zero mass term $m^2$. In the boundary CFT, the scaling dimension of the dual boundary operator $\mathcal{O}$ is computed via the relation
\begin{align}
\label{massRelation}
m^2 = p^{1-\Delta}+p^\Delta-p-1\,.
\end{align}
Solving this equation for real $\Delta$ as a function of $m^2$ yields two solutions $\Delta_+$ and $\Delta_-$ with $\Delta_++\Delta_-=1$. Letting $\Delta = \Delta_+$ denote the larger of the two solutions, the scaling dimension of $\mathcal{O}$ is furnished by $\Delta_+$, while $\Delta_-$ is the dimension of $\phi$, the boundary field obtained from the limiting value of the bulk field.

The free two-point function of $\mathcal{O}$ was determined in \cite{HMST} and also studied in \cite{HJ2024}. In both references, separate answers were given for the two cases that arise depending on whether the two insertion points have equal or distinct norms. But the answers can be combined into a single formula, given, in the notation of the present paper, for two boundary points $x_1$ and $x_2$ in the fundamental $E$, by
\begin{align}
\label{OO}
\left<\mathcal{O}_\Delta(x_1)\mathcal{O}_\Delta(x_2)\right>
=
\frac{|x_1|^\Delta\,|x_2|^\Delta}{|x_1-x_2|^\Delta}
+\frac{1}{p^{m\Delta}-1}
\bigg(\frac{|x_1|^\Delta}{|x_2|^\Delta}+\frac{|x_2|^\Delta}{|x_1|^\Delta}\bigg)\,.
\end{align}
We note that on setting $\Delta=1$, the two point function \eqref{OO} becomes identical to $H(x_1,x_2)$, which of course is not coincidence. Since the exponentiated action $e^{-S}$ for $\phi$ by the AdS/CFT dictionary is equal to the path integral over the dual $\mathcal{O}$ field when sourced with $\phi$, the kernel $H$ for the $\phi$ action equals the inverse kernel for the $\mathcal{O}$ action, ie. the two-point function for $\mathcal{O}$.

In Section~\ref{sec:Green} we identified the the N\'eron-Tate local height function as the Green's function for the differential operator $D$ in the $\phi$ action in the massless limit when $\phi$ is dimensionless. To investigate if the N\'eron-Tate local height function too can be uncovered in $p$-adic thermal CFT correlators, we can consider formally continuing the scaling dimension $\Delta$ beyond the regime of validity for the holographic computation that resulted in \eqref{OO} and study the $\Delta \rightarrow 0$ limit, where $\mathcal{O}$ becomes dimensionless (a limit we also commented briefly on in \cite{HJ2024}). In this limit, the correlator \eqref{OO} has the expansion
\begin{align}
\nonumber
&\hspace{30mm}\left<\mathcal{O}_\Delta(x_1)\mathcal{O}_\Delta(x_2)\right>
=\frac{2}{\Delta\,m\, \log p}+
\\[-2.5mm]  \label{lim1}
\\[-2.5mm] \nonumber
&2\log p\,\bigg(
\hspace{-1mm}-\frac{\log |x_1-x_2|}{\log p}
+\frac{\log(|x_1|\,|x_2|)}{2\log p}
+\frac{\big(\log |x_1|-\log|x_2|\big)^2}{2m(\log p)^2}
+\frac{m}{12}\bigg) + \mathcal{O}(\Delta)\,.
\end{align}
The leading divergent piece is simply a constant. Meanwhile the zero-th order term in $\Delta$ gives a logarithmic correlator that can be checked to be precisely of the form of the N\'eron-Tate height function $h$ in \eqref{h}, including the constant term $\frac{m}{12}$, with the argument of $h$ given by $x_1/x_2$ or equivalently $x_2/x_1$. Hence, we see that the dimensionless limit of the $p$-adic thermal CFT two-point function precisely yields the N\'eron-Tate height function:
\begin{align}
\lim_{\Delta\rightarrow 0} \bigg(\left<\mathcal{O}_\Delta(x_1)\mathcal{O}_\Delta(x_2)\right> -\frac{2}{\Delta\,m\, \log p}\bigg)
=2\log p\, h(x_1/x_2) \,.
\end{align}
\hspace{1mm}
\\

\subsection*{Acknowledgment}

The work of A.~H. is supported by Simons collaboration grant No. 708790. The work of C.~B.~J. is supported by the Korea Institute for Advanced Study (KIAS) Grant PG095901.

\bibliographystyle{amsplain}
%\bibliography{references}
\providecommand{\bysame}{\leavevmode\hbox to3em{\hrulefill}\thinspace}
\providecommand{\MR}{\relax\ifhmode\unskip\space\fi MR }
\providecommand{\MRhref}[2]{%
  \href{http://www.ams.org/mathscinet-getitem?mr=#1}{#2}
}

\end{document}